\documentclass[11pt, oneside]{article}   	% use "amsart" instead of "article" for AMSLaTeX format
\usepackage{lipsum} 
\usepackage{sidecap}
\usepackage{float}
\usepackage{cite}
 \usepackage{cleveref}
\usepackage{color}

\usepackage{lipsum}
\usepackage{mathtools}
\usepackage{amscd}
\usepackage{enumitem}
\usepackage[a4paper, total={5.2in, 9.45in}]{geometry}
\usepackage[usenames,dvipsnames,svgnames,table]{xcolor}
\usepackage{amsmath,amssymb,latexsym}
\usepackage[utf8]{inputenc}
\usepackage[english]{babel}
\usepackage{booktabs}

\newtheorem{thm}{Theorem}[section]
\newtheorem{conj}{Conjecture}[section]

\newtheorem{defi}{Definition}[section]

\newcommand{\PSLZ}{\operatorname{PSL}(2,\mathbb{Z})}
\newcommand{\Cov}{\operatorname{Cov}_0^{Q}}

\newcommand{\parr}[1]{\vspace{0.1cm} \noindent {\bf #1}}

\usepackage{tocloft}
\usepackage{siunitx}
\usepackage{xcolor}
\usepackage{booktabs,colortbl, array}
\usepackage{pgfplotstable}
\definecolor{rulecolor}{RGB}{0,71,171}
\definecolor{tableheadcolor}{gray}{0.92}
 \def\C{\mathbb{C}}

\def\F{\mathcal{F}}

\def\R{\mathbb{R}}

\def\Z{\mathbb{Z}}
\def\H{\mathbb{H}}

\newcommand{\nid}{\noindent}

\newcommand\quotient[2]{
        \mathchoice
            {% \displaystyle
                \text{\raise1ex\hbox{$#1$}\Big/\lower1ex\hbox{$#2$}}%
            }
            {% \textstyle
                #1\,/\,#2
            }
            {% \scriptstyle
                #1\,/\,#2
            }
            {% \scriptscriptstyle  
                #1\,/\,#2
            }
    }

\definecolor{aurometalsaurus}{rgb}{0.43, 0.5, 0.5}
\definecolor{darkjunglegreen}{rgb}{0.1, 0.14, 0.13}
\definecolor{coolblack}{rgb}{0.0, 0.18, 0.39}
\definecolor{cobalt}{rgb}{0.0, 0.28, 0.67}
 
\usepackage{hyperref}
\hypersetup{
     colorlinks   = true,
     citecolor    = cobalt
}

\newcommand{\Addresses}{{% additional braces for segregating \footnotesize
  \bigskip
  \footnotesize

  S.~Bullett, \textsc{School of Mathematical Sciences, Queen Mary University of London,
    London E1 4NS, UK}\par\nopagebreak
  \textit{E-mail address}, S.~Bullett: \texttt{s.r.bullett@qmul.ac.uk}
\medskip

L.~Lomonaco, \textsc{Department of Applied Mathematics, University of S\~ao Paulo,
    Brazil, CEP  05508-900}\par\nopagebreak
  \textit{E-mail address}, L.~Lomonaco: \texttt{lluna@ime.usp.br}

  \medskip

  C.~Siqueira, \textsc{Department of Mathematics, University of S\~ao Paulo,
    Brazil, CEP  05508-900}\par\nopagebreak
  \textit{E-mail address}, C.~Siqueira: \texttt{carloss@ime.usp.br}

  }}

\title{\bf Correspondences in complex dynamics}
\author{Shaun Bullett, Luna Lomonaco and Carlos Siqueira}
\date{}							% Activate to display a given date or no date

\begin{document}

\hypersetup{linkcolor=cobalt}
\maketitle
\begin{quote}
\center{\it In memory of Welington de Melo, whose intellectual honesty inspired us all.}
\end{quote}
\thispagestyle{empty}

%\vspace{0.3cm}
\begin{abstract} 
This paper surveys some recent results concerning the
dynamics of two families of holomorphic correspondences, namely ${\mathcal F}_a:z \to w$ defined by the relation
$$\left( \frac{aw-1}{w-1}  \right)^2 +  \left( \frac{aw-1}{w-1}  \right) \left( \frac{az +1}{z+1}  \right)  + \left( \frac{az+1}{z+1}  \right)^2 =3,$$

\noindent and 
$$\mathbf{f}_c(z)=z^{\beta} +c, \mbox{ where } 1<\beta=p/q \in \mathbb{Q},$$

\noindent which is the correspondence
$\mathbf{f}_c:z \to w$ defined by the relation
$$(w-c)^q=z^p.$$
Both can be regarded as generalizations of the family of quadratic maps $f_c(z)=z^2+c$. We describe dynamical properties for the family $\F_a$
which parallel properties enjoyed by quadratic polynomials, in particular a B\"ottcher map, periodic geodesics and Yoccoz inequality, and we give a detailed account of the very recent theory of holomorphic motions  
for hyperbolic multifunctions in the family ${\bf f}_c$.
\end{abstract}

%\blfootnote{Research partially supported by FAPESP 2016/16012-6.  }
%\pagebreak
\tableofcontents

\thispagestyle{empty}
\newpage

\clearpage
\setcounter{page}{1}
 \section{Introduction}

A holomorphic correspondence on the Riemann sphere is a relation $z\mapsto w$  given implicitly by a polynomial
equation $P(z,w)=0.$ Any rational map is an example of a holomorphic correspondence. Indeed, if $f(z)=p(z)/q(z),$ then $w=f(z)$ iff $P(z,w)=0,$ 
where $P(z,w)= wq(z)-p(z).$ In particular, the family of quadratic polynomials $f_c(z)= z^2+c$ (parametrized by $c\in {\mathbb C}$) can be regarded as an analytic family of holomorphic correspondences.
The grand orbits of any finitely generated Kleinian group can also be regarded as those of a holomorphic correspondence.

This paper is concerned with two families of holomorphic correspondences which generalize quadratic polynomials in different ways. The first is the family ${\mathcal F}_a:z \to w$ defined by

\begin{equation} \label{plk} \left( \frac{aw-1}{w-1}  \right)^2 +  \left( \frac{aw-1}{w-1}  \right) \left( \frac{az +1}{z+1}  \right)  + \left( \frac{az+1}{z+1}  \right)^2 =3, \end{equation}
\noindent where $a\in \mathbb{C}$ and $a\neq 1,$ 
introduced in the early nineties by Bullett and Penrose \cite{Bullett1994}. They proved: 

\begin{thm}\label{quadmat}For every $a$ in the real interval $[4,7],$ the correspondence  $\mathcal{F}_a$ is a mating between some quadratic map $f_c(z)=z^2+ c$
and the modular group $\Gamma=\PSLZ,$ \end{thm}
and conjectured that the connectedness locus for this family is homeomorphic to the Mandelbrot set.\\

The second family is 
\begin{equation} \label{fxs} \mathbf{f}_c(z) = z^{\beta} +c,  \ \ c\in{\mathbb C},
\end{equation}

\noindent where $\beta>1$ is a rational number and   $z^{\beta} = \exp  \frac{1}{q}(\log z^p)$.
If $\beta=p/q$  in lowest terms, then each member of the family \eqref{fxs} of multifunctions is a holomorphic correspondence, defined by the relation $(w-c)^q=z^p.$
Hence $\mathbf{f}_c$ maps every $z\neq 0$ to a set consisting of $q$ points.  If $p$ and $q$ are not relatively prime, we shall
use the notation $z^{p/q}+c$ to express the holomorphic correspondence $(w-c)^q=z^p.$  
Thus  $z^2 +c$ and $z^{4/2} +c$ denote different correspondences. \\

In this paper we describe the dynamics of holomorphic correspondences from various perspectives, exploring the concepts of hyperbolicity 
and holomorphic motions for \eqref{fxs} and describing results concerning a B\"ottcher map, periodic geodesics, and a Yoccoz inequality for the family of matings \eqref{plk}. As we shall see, 
the techniques involved in the two studies are independent, but as we have already noted, both families can be viewed as generalizations of the quadratic family, and our techniques for studying them are motivated by the notions of hyperbolicity, external rays, Yoccoz inequalities  and  local connectivity, which are inextricably related to one another in the study of quadratic polynomials $f_c(z)=z^2+c$. For this reason, it will be convenient to start by recalling some well known facts, techniques and open questions concerning this celebrated family of maps.  Excellent sources for details are the books of Milnor \cite{Milnor} and de Faria and de Melo \cite{FM}. An overview of a century of complex dynamics is presented in the article by Mary Rees \cite{Rees}.

\subsection{Dynamics of quadratic maps}\label{bts}

%\paragraph{Density of hyperbolicity.} 
% Given a family of maps, a central and difficult problem in dynamics consists in showing that most maps
% are structurally stable (or hyperbolic), and provide models for maps in the same component determined by the structurally stable regime. 
% A map is called hyperbolic when all its critical points are attracted to an attracting cycle.

Consider the action of $f_c(z)=z^2+c$ on the Riemann sphere $\widehat{\mathbb{C}}$. For any polynomial of degree $d \geq 2$ acting on $\widehat{\mathbb{C}}$, the point $z=\infty$ is a superattracting fixed point. Let $\mathcal{A}_c$ denote its basin of attraction.
The filled Julia set $K_c= K_{f_c}$ is the set of points with bounded orbit, that is $K_c = \widehat{\mathbb{C}} \setminus \mathcal{A}_c$.
The Julia set $\mathcal{J}_c=\mathcal{J}_{f_c}$ is the common boundary of these regions: $\mathcal{J}_c = \partial K_c = \partial \mathcal{A}_c$.
The \emph{Mandelbrot set} $M$ is the connectedness locus of the  family $f_c(z) = z^2 +c$; in other words, the set of all parameters $c \in \mathbb{C}$ such that $\mathcal{J}_c$ is connected. 

On the basin of attraction $\mathcal{A}_c$, the quadratic polynomial $f_c$ is conformally conjugate to the map $f_0(z)=z^2$ by the so-called \textit{B\"ottcher map} $\varphi_c$ (tangent to the identity at infinity). In the case $\mathcal{J}_c$ (or equivalently $K_c$) is connected, the B\"ottcher map extends to a conformal conjugacy:
$$\varphi_c: \mathbb{C} \setminus K_c \to  \mathbb{C}\setminus \overline{\mathbb{D}}_{1}$$
(An analogue of this map for the family $\F_a$ will appear in Section \ref{lmv}.) 
The \emph{external ray} $R_{\theta}^c\in \mathbb{C} \setminus K_c$ with argument ${\theta} \in \R/\Z$ is the preimage under the B\"ottcher map $\varphi_c $ of the half-line $t e^{2\pi i \theta}\in \mathbb{C}\setminus \overline{\mathbb{D}}_{1}$, with $t \in (1, \infty)$. 
When $$\lim_{t\to 1^+}\varphi_c^{-1} (t e^{2\pi i \theta})=z,$$ we say that $R_{\theta}^c$ lands at $z.$ 
We know that \textit{rational rays land} \cite{DH84, Milnor},  and that \textit{repelling and parabolic periodic points are landing points}
of at least one and at most finitely many rays \cite{Milnor}. By Carath\'eodory's theorem, if $\mathcal{J}_c$ is locally connected, then every external ray lands. We remark that the B\"ottcher map and external rays can also be defined for degree $d$ polynomials, and in this case as well rational rays land and repelling and parabolic periodic points are landing points \cite{Milnor}.
(\textit{Hyperbolic geodesics} play an analogous role for the family $\F_a$ and enjoy similar properties to external rays, see Section \ref{lmv}).\\

Using the B\"ottcher map, Douady and Hubbard constructed a conformal homeomorphism between the complement of the Mandelbrot set and the complement of the closed unit disk:
$$\Phi: \mathbb{C} \setminus M \to  \mathbb{C}\setminus \overline{\mathbb{D}}_{1}$$
$$\,\,\, c \,\,\,\, \rightarrow \,\,\, \varphi_c(c),$$
proving that the Mandelbrot set is compact and connected \cite{DH84}. This isomorphism also allows the definition of 
\emph{ parameter space external rays}: the parameter ray of argument $\theta$ is $\mathcal{R}_{\theta}=\Phi^{-1}(R_{\theta}^0)$. 
If $M$ is locally connected, then every external ray lands. Conjecturally, the Mandelbrot set is locally connected (which we write MLC).
This topological conjecture is crucial in one dimensional complex dynamics, since it has been proved (\cite{DH85}) to imply 
\textit{density of hyperbolicity} for the quadratic family. 
A rational map is called \textit{hyperbolic} when all its critical points are attracted to attracting cycles.
Hyperbolic maps are among the best understood rational maps. Indeed, if the quadratic polynomial $f_c$ is  hyperbolic then (i)  every orbit in the interior of the filled Julia set $K_c$ 
(if non-empty) converges to the finite attracting cycle (which is unique since $f_c$ is quadratic); (ii) every orbit outside $K_c$ converges to $\infty;$ and (iii) $f_c$ is expanding and topologically 
mixing on the Julia set $\mathcal{J}_c =\partial K_c$. 
A major conjecture in holomorphic dynamics is:
 \begin{conj}[Density of hyperbolicity]\label{denhyp}
  The set of hyperbolic
rational maps is open and dense in the space of rational maps $\operatorname{Rat}_d$  of the same degree. 
 \end{conj}
A version of this conjecture dates back to Fatou, and for this reason Conjecture \ref{denhyp} is often known as the \textit{Fatou conjecture}. Note that it concerns density of hyperbolicity, since openness of the set of hyperbolic maps is known.\\

 Strongly related to hyperbolicity is the concept of structural stability. A map $f_a$ is \emph{structurally stable}  if $f_c$ is topologically conjugate to $f_a,$ for every $c$ in an open set containing $a.$ For rational maps on the Riemann sphere \textit{$J$-stability}, which roughly speaking means stability on a neighborhood of the Julia set, is usually considered \cite{Rees}. Ma\~n\'e, Sad and Sullivan \cite{Mane1983} have shown that the set of $J$-structurally stable
rational maps is open and dense in the space of rational maps $\operatorname{Rat}_d$  of the same degree. Since in any family of holomorphic maps
the set of hyperbolic parameters forms an open and dense subset of the $J$-stable parameters, Conjecture \ref{denhyp} is equivalent to the following (see \cite{McMullen94}):

\begin{conj}
 A $J$-stable rational map of degree $d$ is hyperbolic.
\end{conj}

%(\textcolor{red}{where topologically mixing means bla bla})    
%The \emph{Fatou Conjecture} is that among rational maps of any given degree $d>2$, an open dense set are hyperbolic. [\textcolor{red}{Is this %the version we want?}]
% A hyperbolic component is a connected component of the interior of $M$ containing a parameter for which the corresponding map is hyperbolic.
% On the other hand, a \emph{ghost component} (or queer component) is a connected component of the interior of $M$ that is not hyperbolic. 
% Under the hypothesis that $M$ is locally connected, Douady and Hubbard  showed that every ghost component has at least three external arguments on its boundary, which is a cetral fact in their combinatorial argument proving that (\cite{DH85}): 
% \begin{thm}[MLC $\Rightarrow$  Density of Hyp 2]\label{pqr} If  $M$ is locally connected, then hyperbolic maps are open and dense in the quadratic family. 
% \end{thm}

For quadratic polynomials, Conjecture \ref{denhyp} claims that the set of $c$ such that $f_c(z)=z^2+c$ is hyperbolic is an open and dense subset of the complex plane. On the other hand, density of $J$-stability implies that each of the infinitely many components $U$ of $\mathbb{C} \setminus M$ is the
parameterization domain of a holomorphic motion $h_c: \mathcal{J}_a \to \mathcal{J}_c,$ $c\in U$ (holomorphic motions are defined in Section \ref{dfc}), with base point $a\in U$ arbitrarily fixed, and every $h_c$ being a \emph{quasi-conformal conjugacy}. If $U$ is a component of $\mathbb{C}\setminus \partial M$ having one point $a$ for which $f_a(z)= z^2+ a$ is hyperbolic, 
 then $f_c$ is hyperbolic for every $c$ in $U$, and thus in the quadratic setting density of hyperbolicity is equivalent to  
conjecturing that \textit{every component of $\mathbb{C}\setminus \partial M$ is hyperbolic}.  
Note that, since $\mathcal{J}_0=\mathbb{S}^1,$ it follows that $\mathcal{J}_c$ is a \emph{quasicircle} (image of $\mathbb{S}^1$ under a quasiconformal homeomorphism) for every $c$ close to zero (more precisely, for every $c$ in the same hyperbolic component as $c=0$).
(A generalization of this fact for $z^{\beta} +c$ is given by Theorem \ref{bfl}). 

%  Thus, by Theorem \ref{pqr} from Douady and Hubbard, MLC the following is an equivalent form of MLC conjecture:
% 
%  \begin{conj} \label{hgf} Every component of $\mathbb{C}\setminus \partial M$ is hyperbolic. 
% 
% \end{conj}

% 
% 
% The components  of the Mandelbrot set meeting $\mathbb{R}\cup\{\infty\}$ are known to be hyperbolic, but the Fatou Conjecture (Conjecture \ref{hgf}) has proved 
% to be one of the hardest open problems in dynamics for the last hundred years. 

In the late eighties J.-C. Yoccoz made a major contribution towards the MLC conjecture, proving that MLC holds at every point $c \in \partial M$ such that
$f_c$ is not infinitely renormalizable. A key ingredient is what is now known as the \emph{Yoccoz inequality}. It can be shown that if $z$ is a repelling fixed point for a degree $d$ polynomial $P$ with connected filled Julia set, then just finitely many external rays $\gamma_i$, say $q'$, land at $z$. Each $\gamma_i$ is periodic with the same period, and there exists $p'<q'$ such that $P\circ \gamma_i = \gamma_{i+p'}$ for any $i$. The number of cycles of rays landing at $z$ is $m=\operatorname{gcd}(p',q')$, and $\theta= p/q = (p'/m)/(q'/m)$ is called the \emph{combinatorial rotation number} of $P$ at $z$. 

\begin{thm}[Yoccoz-Pommerenke-Levin inequality \cite{Hubbard91, Pommerenke86, Levin91}] \label{yocin} If $z$ is a repelling fixed point of a degree $d$ polynomial $P$ with connected filled Julia set, and $\theta=p/q$  is its combinatorial rotation number in lowest terms, then
\begin{equation}
\frac{\operatorname{Re} \tau}{ |\tau - 2\pi i \theta |^2 }\geq \frac{mq}{2\operatorname{log} d},
\end{equation}
 for some branch $\tau$ of $\operatorname{log} P'(z).$ 
\end{thm}
(A Yoccoz inequality for the family $\F_a$ is developed by the first two authors in \cite{BL17}; see Theorem \ref{xxx}. While the original Yoccoz inequality is proven for degree $d$ polynomials, and so applies to iterates of degree $2$ polynomials and hence to periodic orbits, an inequality of the form presented in Theorem \ref{xxx} has so far only been proved for repelling fixed points.)\\ 

In 1994, C. McMullen made a deep contribution toward MLC, by proving that every component of the interior of the Mandelbrot set meeting the real axis is hyperbolic \cite{McMullen94}.
In the late nineties, M. Lyubich \cite{Lyubich97}, and independently Graczyk and Swiatek \cite{Swiatek98} proved density of hyperbolicity for the real quadratic family. About ten years later Kozlovski, Shen and van Strien proved it for real polynomials of higher  degree, by proving that any real polynomial can be approximated  by hyperbolic real polynomials of the same degree \cite{SSK07}. However, density of hyperbolicity for degree $d$ rational maps on $\widehat{\mathbb C}$ is still open.

\subsection{Dynamics of holomorphic correspondences}

We now outline our main results described in this paper, concerning the families \eqref{plk} and \eqref{fxs}: these involve generalizations of the concepts presented in Section \ref{bts}. {\sl Readers who want to see the  proofs - as Welington always did - } can find those concerning family \eqref{plk} in \cite{BL16} and \cite{BL17}, and those concerning family \eqref{fxs} in \cite{SS17, Carlos, Rigidity, Siq17}.

\paragraph{Part I.} We start with an abstract definition of matings between quadratic maps and $\PSLZ$  (Section \ref{wxz}) with the help of Minkowski's question mark function. This description dates back to 1994, when the first author together with C. Penrose \cite{Bullett1994} started investigating the family $\mathcal{F}_a.$ 
The formal definitions of limit sets and the connectedness locus $\mathcal{C}_{\Gamma}$ for this family are given in Section \ref{gra}.
There we also define a \textit{mating} between the modular group and a map in the parabolic quadratic family
$$\operatorname{Per}_1(1) = \{ P_A(z)= z+1/z+A \, | \, A \in \C \}/(A\sim -A),$$
and present a result which is a significant advance on Theorem \ref{quadmat}, namely that for any $a\in \mathcal{C}_{\Gamma},$ the correspondence $\mathcal{F}_a$ is a mating between $\PSLZ$ a parabolic map in
$\operatorname{Per}_1(1)$ (see Theorem \ref{Mat}, and figures \ref{ppp} and \ref{bbb}). 

We open Section \ref{lmv} by recalling the existence of a B\"ottcher map for the family $\mathcal{F}_a$ when $a \in \mathcal{C}_{\Gamma}$ (see Theorem \ref{efv}), and we then use it to construct periodic geodesics on the regular domain of $\mathcal{F}_a$ (an analogue of periodic external rays). These land (see Theorem \ref{land}), analogously to the rational external rays for the quadratic family of polynomials.

By a quite technical and deep argument \cite{BL17}
it can be shown that when $a$ is in ${\mathcal C}_\Gamma$ every repelling fixed point $z$ of ${\mathcal F}_a$ is
the landing point of exactly one periodic cycle of geodesics. It follows, as for polynomials, 
that $z$ has a well-defined combinatorial rotation number $\theta=p/q$. A geodesic in
the cycle is stabilized by a Sturmian word $W_{p/q}$, in $\alpha$ and $\beta$, of rotation
number $p/q$ (Sturmian words are defined in Section \ref{fp}: $W_{p/q}$ is unique up to cyclic permutation for any given $p/q$). 
%In \cite{BL17} the first two authors prove:

\begin{figure}[H]
\centering
\includegraphics[scale=0.5]{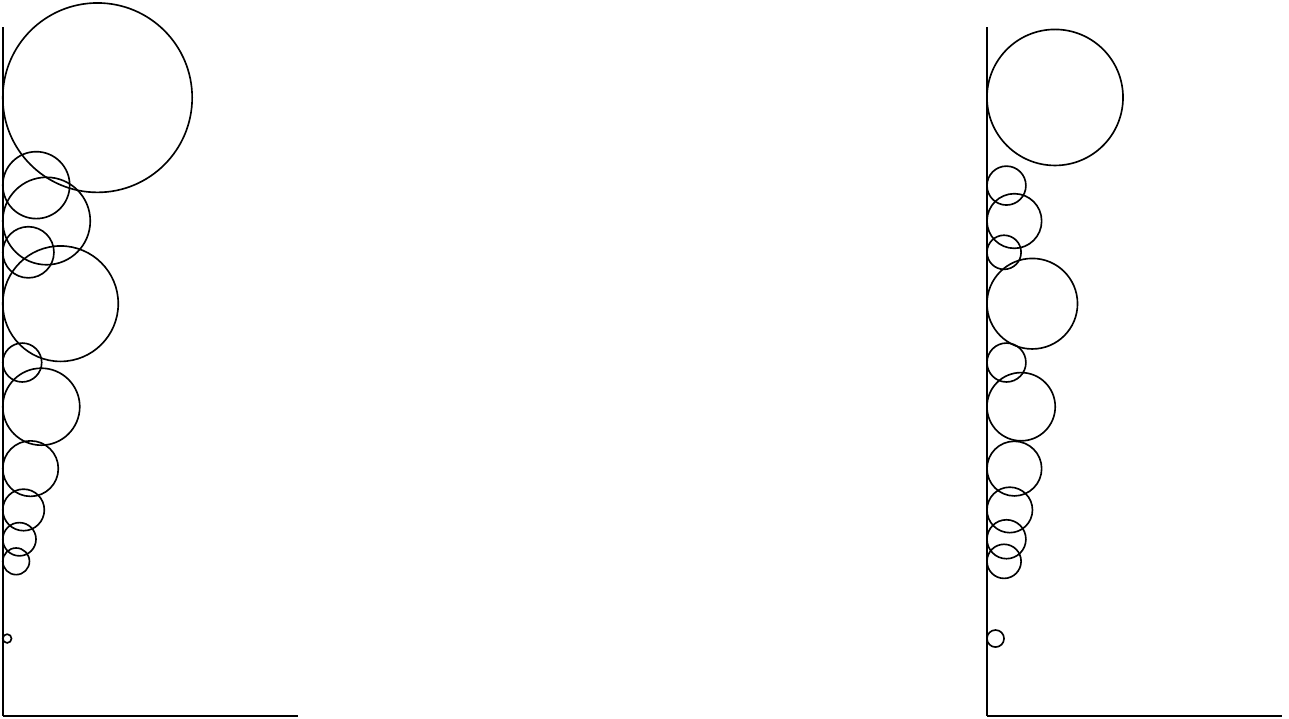}
\caption{Disks in the $\tau$-plane permitted by the Yoccoz inequality: on the left for the matings $\mathcal{F}_a$, and on the right for the classical case of quadratic polynomials. (In each case the disks plotted correspond to all $p/q \in [0, 1/2]$ with $q \leq 8,$ and to $1/16$).} 
\label{Ydiscs}
\end{figure}

%We then associate a Sturmian word $W_{p/q}$ in $\alpha$ and $\beta$ to a periodic geodesic 
%$\gamma$ landing at a repelling fixed point $z$ (Sturmian words are defined in Section \ref{fp}). 
%The rotation number, $p/q$, of this word is equal to the
%combinatorial rotation number of the branch $f_a$ of $\mathcal{F}_a$ which fixes $z$. By a quite technical 
%and deep argument which we will not present here (being even more technical than in the case of 
%polynomials), we obtain that each repelling fixed point is the landing point of exactly one geodesics.
%at least one and at most finitely many periodic geodesics. 
%It follows, as for polynomials, that a repelling fixed point of a correspondence in the family $\F_a$, 
%$a\in {\mathcal C}_\Gamma$, has a well-defined combinatorial rotation number $\theta=p/q$. We can now state:

\begin{thm}[Yoccoz inequality]\label{xxx} Let $a\in \mathcal{C}_{\Gamma}$ and $z$ be a repelling fixed point of $f_a$ whose combinatorial rotation number is $\theta=p/q$ in lowest terms. Then there is a branch $\tau$ of $\operatorname{log} f_a'(z)$ such that 
$$\frac{\operatorname{Re} \tau}{|\tau - 2\pi i \theta |^2}   \geq \frac{q^2}{4p \operatorname{log} (\left \lceil{q/p}\right \rceil 
 +1)},   \ \ \ \textrm{if}   \ \ \theta \leq 1/2;  \ and $$

 $$\frac{\operatorname{Re} \tau}{|\tau - 2\pi i \theta |^2}   \geq \frac{q^2}{4(q-p)\operatorname{log} (\left \lceil{q/(q-p)}\right \rceil 
 +1)},   \ \ \ \textrm{if}   \ \ \theta > 1/2. $$

\end{thm}

The inequalities of both Theorems \ref{yocin} and \ref{xxx} have geometric interpretations as restricting the logarithm of the derivative at a repelling fixed point to a round disk for each $p/q$. 
See Figure \ref{Ydiscs} for illustrations.\\

\begin{figure}[H]
\centering
\includegraphics[scale=0.4]{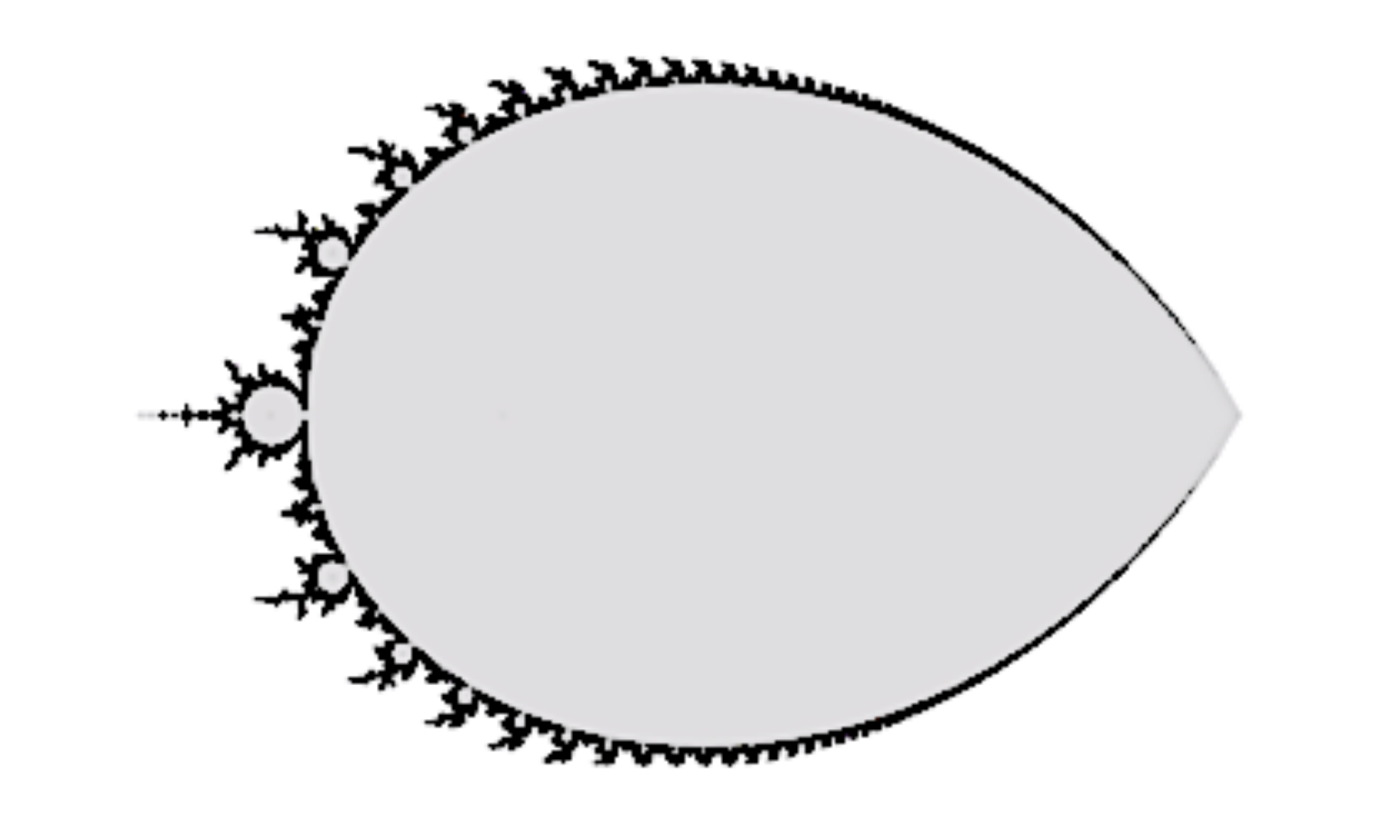}
\caption{ A plot of $M_{\Gamma}=\mathcal{C}_{\Gamma} \cap \{z: |z-4|\leq 3\},$ which is conjecturally homeomorphic to the Mandelbrot set \cite{Bullett1994}.}\label{MGamma}
\end{figure}

Theorem \ref{xxx} provides a key step in the strategy of the first two authors to prove that
the part $M_{\Gamma}=\mathcal{C}_{\Gamma} \cap \{z: |z-4|\leq 3\}$ of the connectedness locus $\mathcal{C}_{\Gamma}$ of the family \eqref{plk} is 
homeomorphic to the connectedness locus $M_1$ of the parabolic family $\{z\mapsto z+ \frac{1}{z} +A:  A\in \mathbb{C}\}/(A\sim -A)$.
With the result announced by Carsten Peterson and Pascale Roesch that $M_1$ is homeomorphic to the Mandelbrot set $M$ \cite{Petersen17}, this will finally prove 
the long-standing conjecture that $M_\Gamma$ (pictured in Figure \ref{MGamma}) is homeomorphic to $M$.

%\begin{thm} The following  three connectedness loci are homemorphic:
%\begin{enumerate} \item the Mandelbrot set $M$ of the quadratic family;
%\item the connectedness locus of the parabolic family $$z\mapsto z+ \frac{1}{z} +A,  \ \ A\in \mathbb{C};$$
%\item  the bounded part $M_{\Gamma}=\mathcal{C}_{\Gamma} \cap \{z: |z-4|\leq 3\}$ of the connectedness locus $\mathcal{C}_{\Gamma}$ of the family \eqref{plk}. 
%\end{enumerate}
%\end{thm}

\paragraph{Part II.} The last Section \ref{tra} describes the dynamics of hyperbolic correspondences in the family \eqref{fxs}.  We start by defining Julia sets (see Figure \ref{fff} for an example). The main subject is the generalization of holomorphic motions, which involves the construction of a solenoid associated to the Julia set of $\mathbf{f}_c(z)=z^{\beta} +c$ (Theorem \ref{bfl}). For parameters $c$ close to zero, the dynamics of $z^{\beta}+c$ on its Julia set $\mathcal{J}_c$ is the projection of a (single-valued) dynamical system $f_c: U \to U$ given by as holomorphic map defined on a subset $U \subset \mathbb{C}^2.$ The maximal invariant set of $f_c$ is a solenoid whose projection is $\mathcal{J}_c.$ The projection of the holomorphic motion in $\mathbb{C}^2$ yields a branched holomorphic motion on the plane, as defined by Lyubich and Dujardin \cite{Lyubich2015} for polynomial automorphisms of $\mathbb{C}^2.$ Branched holomorphic motions  are described in greater generality for the family \eqref{fxs}
 in  \cite{SS17}.

The advantage of the solenoid construction is that it makes possible to apply certain techniques of Thermodynamic  Formalism to the family of maps $f_c:U \to U$ and use them to estimate the Hausdorff dimension of $\mathcal{J}_c.$ For example,  

 \begin{thm}[Hausdorff dimension]\label{tyu} If $q^2 < p$ then for every $c$ sufficiently close to zero,
$$\dim_{H} \mathcal{J}_c < 2, $$

\noindent where $\dim_H$ denotes the Hausdorff dimension of $\mathcal{J}_c.$

\end{thm}

\begin{figure}[H] 
\centering
\includegraphics[scale=0.43]{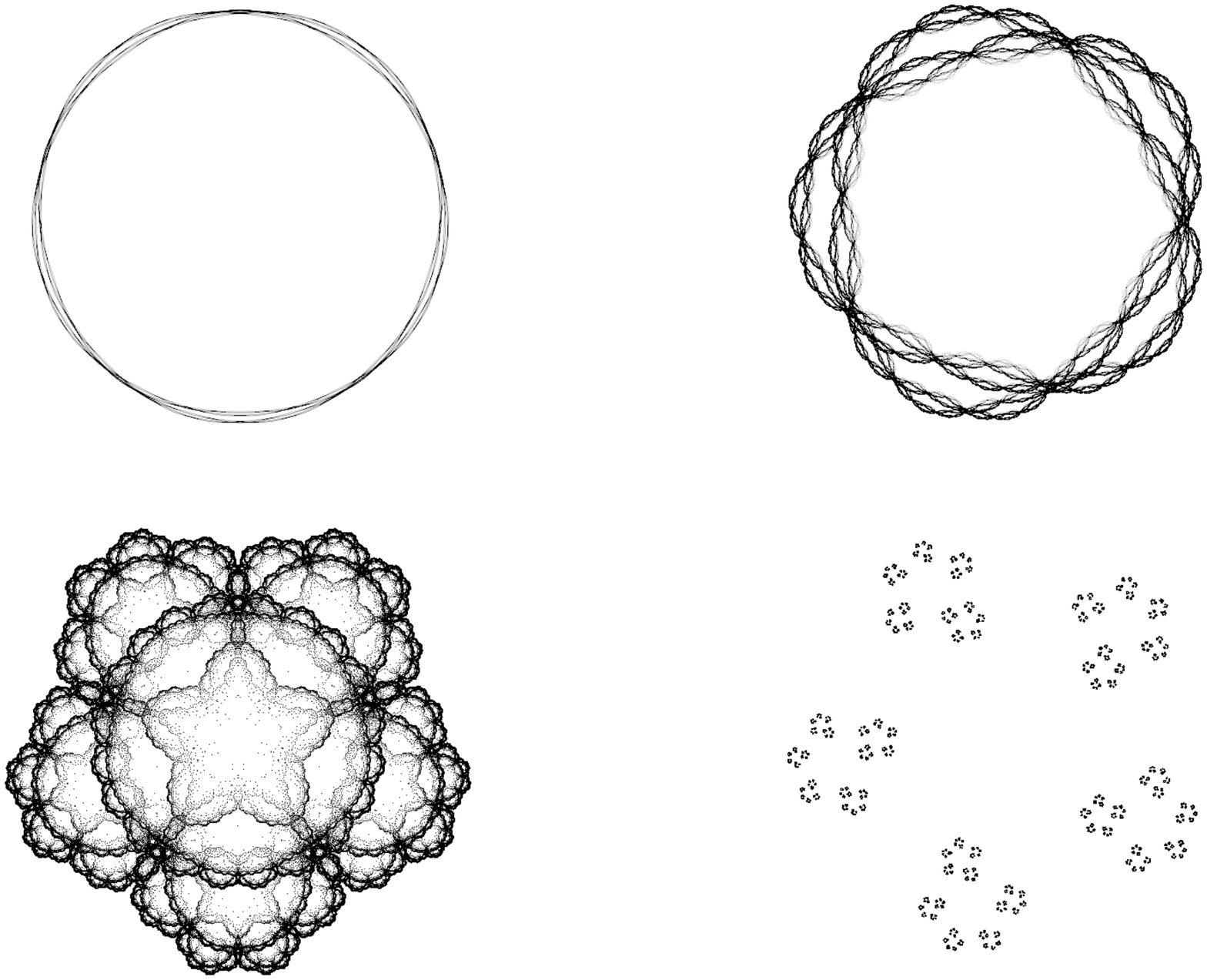} 
\caption{  Julia sets of $z^{\beta} +c,$ where $\beta=5/2$ is fixed.  The values  of $c$ are, respectively,  $0.05, (1+i)/5,\  0.7$ and $2 +i,$ read from upper-left to bottom-right.     It starts as a circle at the singularity $c=0,$ but the first figure reveals that $\mathcal{J}_c$ is the \emph{shadow} of a solenoid for every $c$ close to zero with $c\neq 0.$ As we perform the branched motion, more bifurcations are added to $\mathcal{J}_c.$ The complexity increases up to a certain moment (third to fourth steps) when the process reverses and  $\mathcal{J}_c$ becomes a Cantor dust.  
  The first three are connected and the fourth  is a Cantor set.  In this family,  $\mathcal{J}_c$ is a Cantor set for  $|c|$ sufficiently large. 
}
\label{bcx}
\end{figure}

In the family of Figure \ref{bcx} we have $p=5$ and $q=2.$ Since $2^2 < 5,$ it follows that $\mathcal{J}_c$ is the projection of a solenoid having zero Lebesgue measure.  The assumption $q^2 < p$ may not be sharp. The essential idea is that $\dim_H\mathcal{J}_c \to 2$ as $\beta \to 1,$ which is supported by many experiments.   

\paragraph{Notation and terminology.}  \begin{enumerate} \item Holomorphic correspondences are denoted by $\mathcal{F}, \mathcal{G}, \ldots$  in the context of matings, or by  $\mathbf{f}, \mathbf{g}, \ldots$ when studying hyperbolic 
%correspondences. 
multifunctions.

\item  %We shall use \emph{iff} instead of \emph{if, and only if.} 
By the term \emph{multifunction} we mean any multivalued map. Every multifunction maps points to subsets. 

 \item $\mathbb{S}^1= \{z\in \mathbb{C}: |z|=1\},$ $\widehat{\mathbb C} = \mathbb{C}\cup \{\infty\},$ $\mathbb{H}= \{z=x+iy \in \mathbb{C}:y>0\}, $ and $f^n=  \underbrace{f\circ \cdots\circ f}_{n}.$

 \item $\Gamma=\PSLZ$ is the modular group consisting of all M\"obius transformations $$z\mapsto  \frac{az+b}{cz +d}, $$
 \noindent where $ad-bc=1$ and $a,b,c,d \in \mathbb{Z}.$ The operation is the standard  composition  $\circ$. The \emph{generators} of the modular group 
that we shall use are 
 the maps $$\alpha(z)=z+1 \mbox{ and } \beta(z)=\frac{z}{z+1}. $$ 
 
 Consider
\begin{equation}\label{scv}  P(z,w) =  (w-(z+1))(w(z+1) -z)=0.\end{equation} 

The grand orbits of $\PSLZ$ on ${\mathbb H}$ are identical to those of  the holomorphic correspondence $\mathcal{H}: \widehat{\mathbb C} \rightarrow \widehat{\mathbb C}$ determined by $P(z,w)=0.$ 
\end{enumerate}

   \parr{Acknowledgments.} 
The authors would like to thank the Funda\c{c}\~ao de amparo a pesquisa do estado de S\~ao Paulo, which has supported the first two authors by the grant FAPESP 2016/50431-6, and the third by FAPESP 2016/16012-6. C.S. is very grateful to Edson de Faria for the hospitality at IME-USP,  to Sylvain Bonnot for suggesting the investigation of some interesting problems relating holomorphic correspondences to automorphisms of $\mathbb{C}^2,$ and to Daniel Smania for many discussions and key ideas on the dynamics of hyperbolic correspondences, specially those concerning Gibbs states and Hausdorff dimension. 

 L.L. and C.S. would like to express their sincere gratitude to the scientific committee and organizers of the conference \emph{New trends in one-dimensional dynamics}, on the occasion of Welington de Melo's seventieth birthday, specially to Pablo Guarino and Maria J. Pacifico (significant content of this paper has been previously announced in this conference).

\section{Mating quadratic maps with $\PSLZ$}

% Thurston's characterisation of rational maps among critically finite branched coverings of the sphere provides a criterion to decide whether two hyperbolic quadratic polynomials are \emph{matable}.    

Recall that in the case of hyperbolic quadratic polynomials $f_c(z)=z^2+c$, the \emph{topological mating} between $f_c$ and $f_{c'}$ is the map $$g: \frac{K_c \cup K_{c'}}{\sim}   \to  \frac{K_c \cup K_{c'}}{\sim}  $$ induced by $f_c$ and $f_{c'}$ on the quotient space, where $\sim$ is the smallest closed relation such that $\varphi_c(z) \sim \varphi_{c'}(\overline{z})$, for every $z\in \mathbb{S}^1$ ($\varphi_c$ is the boundary extension of the B\"ottcher coordinate and $K_c$ is  a copy of the filled Julia set). The two maps are \emph{matable} if  the quotient space 
is a sphere, and $g$ can be realized as a rational map. By applying Thurston's characterization of rational maps among critically finite branched coverings of the sphere, Tan Lei (\cite{Tan}) and Mary Rees (\cite{MR}) proved that two quadratic polynomials $f_c,f_{c'}$ with periodic critical points are matable if and only if $c$ and $c'$ do not belong to
complex conjugate limbs of the Mandelbrot set.

Matings can also be constructed between Fuchsian groups: by applying the Bers Simultaneous Uniformization Theorem certain Fuchsian groups can be mated with (abstractly isomorphic) Fuchsian groups to yield \emph{quasifuchsian} Kleinian groups. (See \cite{Bullett2010} for a discussion of matings in various contexts in conformal dynamics.) What is a surprise when first encountered is that certain Fuchsian groups can be mated  with polynomial maps (see Section \ref{wxz}). This is achieved in a larger category of conformal dynamical systems, containing both rational maps and finitely generated Kleinian groups, the category of holomorphic correspondences on the Riemann sphere. These are multifunctions  $\mathcal{F}: \widehat{\mathbb C} \to \widehat{\mathbb C},$ for which there is a polynomial $P(z,w)$ in two complex variables such that  $\mathcal{F}(z) = \left\{w\in \widehat{\mathbb C}: P(z,w)=0 \right\}$. 
%for every $z\in \widehat{\mathbb C}$. 
%Every rational map or finitely generated Kleinian group can be presented as a correspondence. 

\subsection{Mating quadratic polynomials with $\PSLZ$} \label{wxz}
%Correspondences are objects which generalise both polynomials and finitely generated Kleinian groups 

Examples of matings between quadratic polynomials and 
the modular group were discovered by the first author and 
Christopher Penrose in the early '90s. To understand their existence
we first consider how one can construct an abstract
(topological) model (see also \cite{Bullett1994} and \cite{BodilNuria} for more details).

\paragraph{Topogical mating: Minkowski's question mark function.}

Let $$h:\hat{\mathbb R}_{\ge 0} \to [0,1]$$
denote the homeomorphism which sends 
$x \in {\mathbb R}$ represented by the
continued fraction 
$$ [x_0; x_1, x_2, \ldots] = x_0+\cfrac{1}{x_1+\cfrac{1}{x_2+\cfrac{1}{x_3+\dots}}}$$
to the binary number 
$$ h(x)=0.\underbrace{1\ldots1}_{x_0}\underbrace{0\ldots0}_{x_1}\underbrace{1\ldots1}_{x_2}\ldots$$

This is a version of Minkowski's question mark function \cite{Mink04}.
It conjugates
the pair of maps $\alpha:x \to x+1$, $\beta:x \to x/(x+1)$ to the pair of maps
$t \to t/2$, $t \to (t+1)/2$ (the inverse binary shift). 

If the Julia set $J(f_c)$ of $f_c: z \to z^2+c$ is connected and
locally connected then the B\"ottcher map $\varphi_c:\widehat{\mathbb
C}\setminus {\mathbb D} \to \widehat{\mathbb C} \setminus K(f_c)$ extends to a
continuous surjection $S^1 \to J(f_c)$, which semi-conjugates the
map $z \to z^2$ on $S^1$ (the binary shift) to the map $f_c$ on
$J(f_c)$. We deduce that we may use the homeomorphism $h$
described above to glue the action of $f_c^{-1}$ on $J(f_c)$ to
that of $\alpha$, $\beta$ on $\hat{\mathbb R}_{\ge
0}/_{\{0 \sim \infty\}}$. Equally well we can glue the action of $f_c^{-1}$ on
$J(f_c)$ to that of $\alpha^{-1}$, $\beta^{-1}$ on $\hat{\mathbb
R}_{\le 0}/_{\{0 \sim -\infty\}}$.

\medskip We now take two copies $K_-$ and $K_+$ of the filled Julia set $K_c$ of $f_c$ and glue
them together at the boundary point of external angle $0$
to form a space $K_- \vee K_+$. Each point $z\in K_c$ has a
corresponding $z'$ defined by $f_c(z')=f_c(z)$. Consider the $(2:2)$
correspondence defined on $K_- \vee K_+$ by sending

\smallskip
$\bullet$ $z \in K_-$ to $f_c(z) \in K_-$ and to $z' \in K_+$;

$\bullet$ $z \in K_+$ to $f_c^{-1}(z) \in K_+$.

\smallskip
It is an elementary exercise to check that this correspondence on
$K_- \vee K_+$ can be glued to
%$\mathcal{H}_{|\mathbb{H}}$ (
the correspondence defined by $\alpha$ and $\beta$ 
on the complex upper half-plane using the
homeomorphisms $\hat{\mathbb R}_{\ge 0}/\{0 \sim \infty \} \to
\partial K_-$ and $\hat{\mathbb R}_{\le 0}/\{0 \sim -\infty\}
\to \partial K_+$ defined above. Thus we have a topological mating
between the action of the modular group on the upper half-plane
and our $(2:2)$ correspondence on $K_- \vee K_+$.

\paragraph{Holomorphic mating.}
Reassured by the existence of this topological construction, we define a (holomorphic) mating between a  
quadratic polynomial $f_c,\,\,c \in M$ and 
$\Gamma=PSL(2,\Z)$ to be a $(2:2)$ holomorphic correspondence $\F$ such that:
 \begin{enumerate}
  \item there exists a completely invariant open simply-connected region $\Omega$ and a conformal bijection 
 $\varphi: \Omega \rightarrow \H$ conjugating $\F|_{\Omega}$ to $\alpha|_\H$ and $\beta|_\H$;
 \item $\widehat \C \setminus \Omega = \Lambda= \Lambda_- \cup \Lambda_+$, where $\Lambda_- \cap \Lambda_+=\{P\}$ (a single point)
 and there exist homeomorphisms $\phi_{\pm}: \Lambda_{\pm} \rightarrow K_c$ conjugating respectively
 $\F|_{\Lambda_-}$ to $f_c|_{K_c}$ and  $\F|_{\Lambda_+}$ to $f^{-1}_c|_{K_c}$
 \end{enumerate}
% Let $\mathcal{D}$ denote the set $\{a:|a-4|<3\}\cup\{7\}$ and $M_\Gamma$ the intersection of $\mathcal{D}$ with the 
% connectedness locus of $\Lambda_a=\Lambda_{a,-}\cup\Lambda_{a,+}$ 
% (the set of values of the parameter $a$ for which $\Lambda_a$ 
% is connected). 
In 1994 the first author and C. Penrose proved that  for all parameters $a$ in the real interval $[4,7],$  the correspondence $\F_a$ is a mating between a quadratic polynomial
 $f_c(z)=z^2+c$, $c\in [-2,+1/4]\subset \R$ and the modular group $\Gamma=PSL(2,\Z)$ (see \cite{Bullett1994}).

\subsection{The regular and limit sets of ${\mathcal F}_a$} \label{gra}

Consider the family of holomorphic correspondences $\mathcal{F}_a: \widehat{\mathbb C} \to \widehat{\mathbb C},$
defined by the polynomial equation \eqref{plk}. 
The change of coordinate $\phi_a:\widehat{\mathbb C}  \to \widehat{\mathbb C}$ given by  $$\phi_a(z)= \frac{az+1}{z+1}$$ conjugates $\mathcal{F}_a$ to   the correspondence \begin{equation}\label{pls} J \circ \operatorname{Cov}_0^{Q}, \end{equation} where $J$   is the (unique) conformal involution 
fixing $1$ and $a$, and $Cov_0^Q$ is the  deleted covering correspondence of the function $Q(z)=z^3$, that is to say, the correspondence defined by
the relation
$$\frac{Q(w)-Q(z)}{w-z}=0, \mbox{  i.e.  }z^2+ zw+ w^2=3.$$  
So $\mathcal{F}_a$ and $J \circ \operatorname{Cov}_0^{Q}$ are the same correspondence in different coordinates, and in that sense we write $\mathcal{F}_a= J \circ \operatorname{Cov}_0^{Q}.$

\begin{figure}[H]
\centering
\includegraphics[scale=0.4]{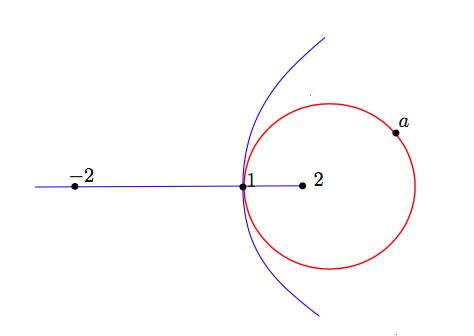}
\caption{Standard fundamental domains for $\Cov$ and $J.$ The curve in blue is $\Cov((-\infty, -2]).$ The region to the right of this curve is a fundamental domain $\Delta_{Cov}$ of $\Cov.$ The unbounded region determined by the red circle is a fundamental domain $\Delta_J$ of the involution $J.$ The parabolic fixed point $P$ is the point $1$.}\label{aaa}
\label{H}
\end{figure}

By a  \emph{fundamental domain} for $\Cov$ (respectively $J$) we mean any maximal open set $U$ which is disjoint from $\Cov(U)$ (respectively $J(U)$).  
We require our fundamental domains to be simply-connected and bounded by Jordan curves (see Figure \ref{H}).

\paragraph{Klein combination locus.}Let $P=1$ denote the common fixed point of $\Cov$ and $J.$ The point $P$ is a \textit{parabolic fixed point}.  The \emph{Klein combination locus} $\mathcal{K}$ is the subset of $\mathbb{C}$ consisting of all $a$ for which there are fundamental domains $\Delta_{Cov}$ and $\Delta_{J}$ of $\Cov$ and $J$, respectively,  such that $$\Delta_{Cov} \cup \Delta_J= \widehat{\mathbb C}\setminus \{P\}.$$
We call such a pair of fundamental domains a \emph{Klein Combination pair}.\\

In \cite{BL16} we show that $\{a\in \mathbb{C}: |a-4|\leq 3,  a\neq 1 \} \subset \mathcal{K}$, and that when
 $a$ is in the interior of this disk the \emph{standard} fundamental domains (see figure \ref{aaa}) are a Klein combination pair. More generally we prove that
for every $a\in \mathcal{K},$ we can always choose a Klein combination pair whose boundaries $\partial \Delta_{Cov}$ and $\partial \Delta_J$ are transversal to the attracting-repelling axis at $P$.

\begin{figure}
  \centering
  \begin{minipage}{.45\textwidth}
  %\centering
  \includegraphics[height= 4.5cm]{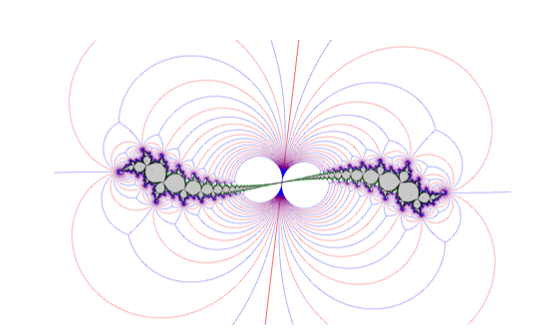}
  \caption{ A connected limit set for $\mathcal{F}_a,$ where $a=4.56 + 0.42i.$}\label{ppp}
  \end{minipage} 
  \hspace{1cm}
  \begin{minipage}{.45\textwidth}
  \centering
  \includegraphics[height= 4.5cm]{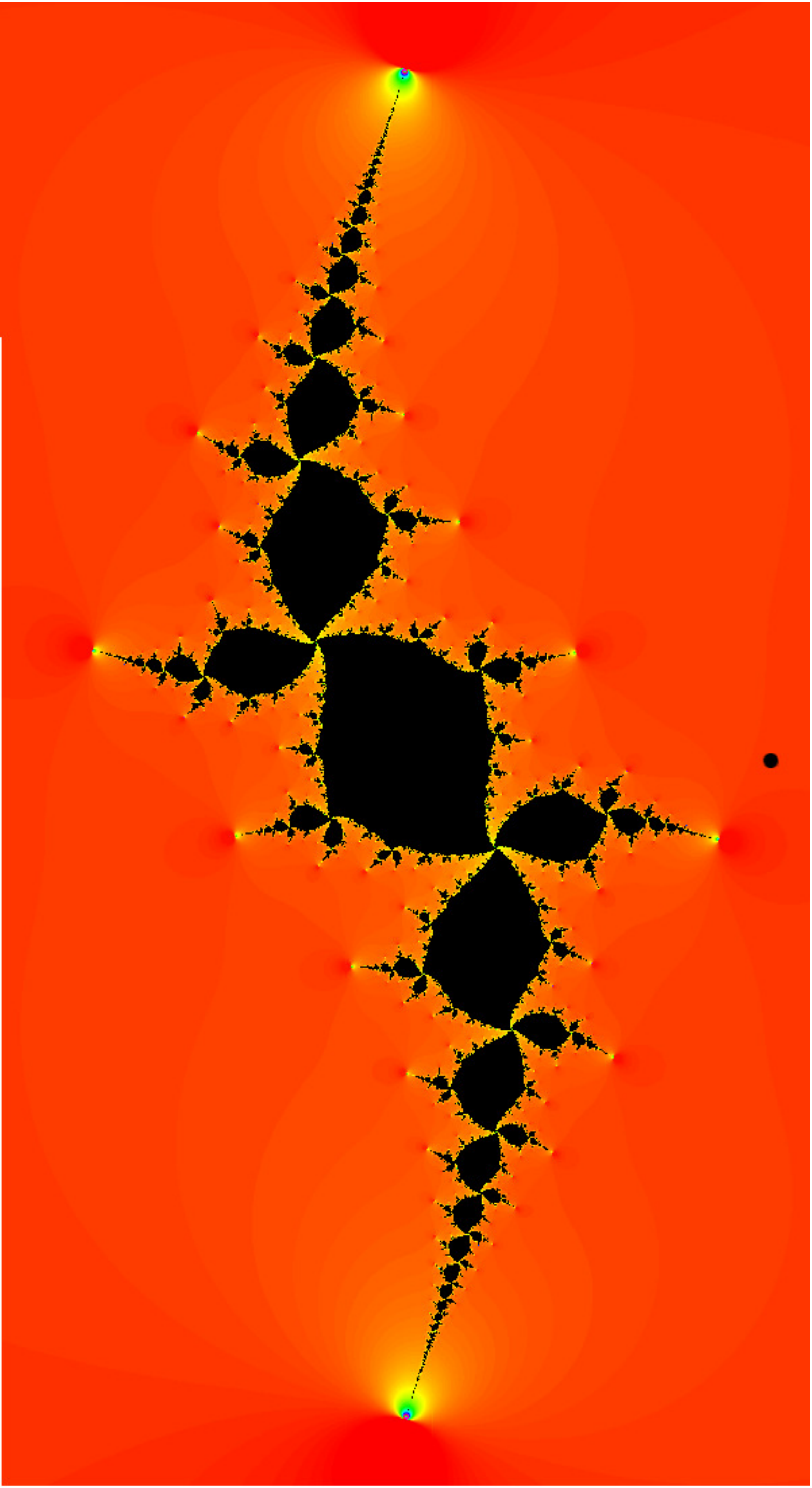}
  \caption{ Julia set of the hybrid equivalent member of $\operatorname{Per}_1(1)$.}\label{bbb}
  \end{minipage}
\end{figure}

Now suppose $a\in \mathcal{K}$ and let $\Delta_{Cov}$ and $\Delta_J$ be a corresponding pair of fundamental domains of $\Cov$ and $J$ such that $\partial \Delta_{Cov}$ and $\partial \Delta_J$ are transversal to the attracting-repelling axis at $P.$ It follows that $P\in \mathcal{F}^{n}_a(\overline{\Delta}),$ for every $n,$ and $\mathcal{F}_a(\Delta)$ is compactly contained in $\Delta \cup \{P\}.$ By definition,\begin{equation}
\Lambda_{a,+} = \bigcap_{n=1}^{\infty} \mathcal{F}_a^n(\overline{\Delta})
\end{equation}
where $\mathcal{F}_a = J \circ \operatorname{Cov}_0^{Q},$
\noindent is the \emph{forward limit set}  of $\mathcal{F}_a.$   Similarly, since $\Delta_{Cov}$ is forward invariant, the complement of $\Delta_{Cov}$ is invariant under $\mathcal{F}^{-1}_a$ and 
\begin{equation}
\Lambda_{a,-}=\bigcap_{n=1}^{\infty} \mathcal{F}_a^{-n} (\widehat{\mathbb C}\setminus \Delta_{Cov})
\end{equation}

\noindent is the \emph{backward limit set} of $\mathcal{F}_a.$  The sets $\Lambda_{a,-} $ and $\Lambda_{a,+}$ have only one point in common, the point $P.$ Their union, $\Lambda_a,$ is the \emph{limit set} of $\mathcal{F}_a.$  An example of a plot of a limit set of $\mathcal{F}_a$ is displayed in Figure \ref{ppp}. (In this plot we use the original coordinate system of \eqref{plk}, so $P=0$ and $J$ is the involution $z \leftrightarrow -z$.) \\

 We have $\mathcal{F}_a^{-1}(\Lambda_{a,-}) =  \Lambda_{a,-},$ and the restriction of $\mathcal{F}_a$ to this set  is a (2\,:1) single-valued holomorphic map denoted by $f_a.$
The involution $J$ maps $\Lambda_{a,-} $ onto $\Lambda_{a,+}$ and determines a conjugacy from  $f_a$    to $$\mathcal{F}_a^{-1}: \Lambda_{a,+} \to \Lambda_{a,+}.$$ 
The \emph{regular domain} of $\mathcal{F}_a$ is $\Omega_a = \widehat{\mathbb C}\setminus \Lambda_a.$ This set is completely invariant under $\mathcal{F}_a$ (forward and backwards). By the Klein Combination Theorem it can be shown that if $\Omega_a$ contains no critical points it is tiled by copies of the intersection 
of any pair of Klein combination  domains, \cite{Bullett00}.

%The  corresponding \emph{limit set} $\varphi^{-1}(\Lambda_a)$ of $\mathcal{F}_a$ is plotted in figure \ref{ppp}, where $\varphi$ is the conjugacy from %$\mathcal{F}_a$ to $\mathcal{G}_a.$ 

%     \begin{figure}[H] 
% \centering
% \includegraphics[scale=0.5]{Limitquad.png} 
% \caption{  A connected limit set for $\mathcal{F}_a,$ where $a=4.56 + 0.42i.$}\label{ppp}
% \end{figure}

\paragraph{Connectedness locus.}
 The \emph{connectedness locus} $\mathcal{C}_{\Gamma}$ of the family $\mathcal{F}_a$ is the subset of $\mathcal{K}$ consisting of all $a$ such that the limit set $\Lambda_a$ is connected. When $a\in \mathcal{C}_{\Gamma}$, the regular domain $\Omega_a$ contains no critical points, and moreover is simply connected.

Bullett and Penrose \cite{Bullett1994} conjectured that for every $a \in \mathcal{C}_{\Gamma},$ the correspondence $\mathcal{F}_a$ is a mating  between some quadratic map $f_c(z)=z^2 +c$ and the modular group $\PSLZ.$ More recently, this conjecture was settled affirmatively by Bullett and Lomonaco \cite{BL16}, provided the quadratic family is replaced by a quadratic family of \textit{parabolic} maps (see figures \ref{ppp} and \ref{bbb}).
 % Let $\mathcal{D}$ denote the set $\{a:|a-4|<3\}\cup\{7\}$ and $M_\Gamma$ the intersection of $\mathcal{D}$ with the 
% connectedness locus of $\Lambda_a=\Lambda_{a,-}\cup\Lambda_{a,+}$ 
% (the set of values of the parameter $a$ for which $\Lambda_a$ 
% is connected). 
% In 1994 the first author and C. Penrose proved that  for all parameters in $M_\Gamma \cap \R$, the correspondence $\F_a$ is a mating between a quadratic polynomial
% $f_c(z)=z^2+c$, $c\in [-2,+1/4]\subset \R$ and the modular group $\Gamma=PSL(2,\Z)$ (see \cite{Bullett1994}).
\subsection{Mating parabolic maps with $\PSLZ.$}
 
 %\paragraph{Mating parabolic maps with $\PSLZ.$} 
The family $\operatorname{Per}_1(1)$ consists of quadratic rational maps of the form $P_A(z) = 1+ 1/z + A,$ where $A\in \mathbb{C}$. The maps in $\operatorname{Per}_1(1)$ all have a persistent parabolic fixed point at $\infty$ and critical points at $\pm 1$. The connectedness locus for the family $\operatorname{Per}_1(1)$ is the \textit{parabolic Mandelbrot set} $M_1$, which has been proved to be homeomorphic to the Mandelbrot set by C. Petersen and P. Roesch (\cite{Petersen17}).
We say that $\mathcal{F}_a$ is a \emph{mating} between $P_A$ and $\PSLZ$ if: 
\begin{enumerate}
\item  on the completely invariant open simply-connected region $\Omega_a$ there exists a conformal bijection 
 $\varphi_a: \Omega_a \rightarrow \H$ conjugating $\mathcal{F}_a: \Omega_a \to \Omega_a$ to $\alpha|_\H$ and $\beta|_\H$; and
\item the (2\,:1) branch of $\mathcal{F}_a$ which fixes $\Lambda_{a,-}$ (given by the holomorphic map $f_a$) is hybrid equivalent  to $P_A$ on the backward limit set $\Lambda_{a,-}$.

\end{enumerate}

In \cite{BL16}, using the theory of parabolic-like maps developed by the second author (see \cite{L}), the first two authors proved  the following (see figures \ref{ppp} and \ref{bbb}):

\begin{thm}\label{Mat}  For every $a\in \mathcal{C}_{\Gamma},$ the correspondence $\mathcal{F}_a$ is mating between a parabolic map in $\operatorname{Per}_1(1)$ and $\PSLZ.$ 

\end{thm}

 The following conjecture has been open for at least 20 years \cite{Bullett1994}:

\begin{conj}\label{MM1} The Mandelbrot set is homeomorphic to $M_{\Gamma}$.   
\end{conj} 

The first two authors have developed a detailed strategy for proving that $M_{\Gamma}$ is homeomorphic to $M_1$.
This, together with the proof by Petersen and Roesch that $M_1$ is homeomorphic to $M$, would finally prove Conjecture \ref{MM1}. 
A key step in the strategy to prove that $M_{\Gamma}$ is homeomorphic to $M_1$ makes use of a Yoccoz inequality for matings, which we prove using a generalization of the technique of external rays (the subject of the next section).

\subsection{Periodic geodesics}  \label{lmv}

%We now generalise the Yoccoz inequality for the family of matings \eqref{plk}. 

\paragraph{B\"ottcher coordinates.}   Consider the holomorphic correspondence $\mathcal{H}$ on the upper half-plane 
obtained from the generators $\alpha(z)=z+ 1$ and $\beta(z)=z/(z+1)$ of $\PSLZ$, i.e. defined by the polynomial equation \eqref{scv}.
As part of the proof of Theorem \ref{Mat} it is shown in \cite{BL16} that:

%Recall that $\Omega_a$ is the regular domain of $\mathcal{F}_a;$  $\mathcal{C}_{\Gamma}$ is the connectedness locus of the family $\mathcal{F}_a;$ and %$\mathbb{H}$ is the hyperbolic plane ($z=x+iy$ for $y> 0$). 

\begin{thm}[B\"ottcher map] \label{efv}If $a\in \mathcal{C}_{\Gamma},$ there is a unique conformal homemorphism $\varphi_a: \Omega_a \to  \mathbb{H}$ 
%tangent to the identity at $\infty$ 
such that $$ \mathcal{H}\circ \varphi_a = \varphi_a \circ \mathcal{F}_a. $$

\end{thm}

By the Schwarz lemma, the B\"ottcher map is an isometry with respect to the hyperbolic metric, and maps geodesics to geodesics. Geodesics in $\Omega_a$, or equivalently in ${\mathbb H}$, play a role for the correspondences ${\mathcal F}_a$ analogous to the role played by external rays for quadratic polynomials $f_c$.

\paragraph{Periodic geodesics land.} By a finite word in $\alpha$ and $\beta$ we mean any M\"obius transformation
  $$W=g_1g_2 \cdots g_n:=g_1 \circ \cdots \circ g_n, $$ where $g_i \in \{\alpha, \beta\}.$  We can compose words in the obvious way, and also consider infinite sequences  $(g_i)_1^{\infty}$  and bi-infinite sequences $(g_i)_{-\infty}^\infty$. 

  A geodesic $\gamma$ in the hyperbolic plane is said to be \emph{periodic} if $W\circ \gamma= \gamma$ for some finite word $W$. (Note that $W$ must include both the letters $\alpha$ and $\beta$, since these being parabolic transformations of ${\mathbb H}$ there are no geodesics invariant under either). 

Since $\mathbb{H}$ is geodesically complete, $\gamma$ is a curve $\mathbb{R} \to \mathbb{H},$ and the limits $$ \gamma(-\infty):=\lim_{t\to -\infty} \gamma(t) , \ \ \gamma(\infty):= \lim_{t\to \infty} \gamma(t)$$ are by definition the \emph{landing points} of $\gamma$. Every periodic geodesic lands on the hyperbolic plane, and the landing points are in $\mathbb{R}\cup \{\infty\}$.  
  
  If $a\in \mathcal{C}_{\Gamma},$  the regular domain is a hyperbolic Riemann surface, that is, it has a unique complete metric of constant curvature $-1$ determining its geometry. A geodesic $\hat\gamma$ in $\Omega_a$ is periodic  if $\varphi_a\circ \hat\gamma$ is a periodic geodesic of $\mathbb{H}.$ 
  
  We say that a periodic geodesic $\hat\gamma:\mathbb{R} \to \Omega_a$  \emph{lands} if the limits $\hat\gamma(\infty)$  and $\hat\gamma(-\infty)$ exist. They are the \emph{right and left landing points,} respectively. 
  
  \begin{thm}\label{land} If $a\in \mathcal{C}_{\Gamma},$ then every periodic geodesic lands. The left landing point belongs to $\Lambda_{a, -}$ and the right landing point is in $\Lambda_{a,+}.$
   \end{thm}

As a corollary, the B\"ottcher map extends to all landing points of periodic geodesics. Indeed it extends to all landing points of preperiodic geodesics, and moreover these correspond under $\varphi_a$ to the set of all quadratic irrationals in ${\mathbb R}$ (the set of real numbers with preperiodic continued fraction expansions).

\subsection{Repelling fixed points, and Sturmian sequences}\label{fp}
The following result is again analogous to a result for quadratic polynomials, but the proof is quite technical and deep (even more so than in the case of polynomials, which is already difficult, see \cite{BL17}), and at present we only have a proof for repelling fixed points, whereas for polynomials it is known for repelling and parabolic cycles:
\begin{thm}\label{lp}
A repelling fixed point in $\Lambda_-({\mathcal F}_a)$ of a correspondence ${\mathcal F}_a$ with $a\in {\mathcal C}_\Gamma$ is the landing point of exactly one periodic cycle of geodesics.
\end{thm}
This theorem has the consequence that to a repelling fixed point $z\in \Lambda_-$ of a correspondence $\F_a$ with $a \in {\mathcal C}_\Gamma$ we can associate a periodic geodesic $\hat\gamma$ which lands there, and a finite word
$W$ in $\alpha$ and $\beta$ which fixes $\varphi_a\circ \hat\gamma$. Letting $f_a$ denote the (locally defined) branch of $\F_a$ which fixes $z$, we deduce that since $f_a$ is 
locally a homeomorphism the cyclic order of the images of $\hat\gamma$ around $z$ is preserved by $f_a$. Thus $f_a$ has a well-defined {\emph{combinatorial rotation 
number} around $z$, and this number is rational.

\paragraph{Sturmian sequences.} %NEW VERSION

Recall that a sequence $(s_i) \in \{0,1\}^{\mathbb{N}}$ is \emph{Sturmian} if, for every $n,$ the number of $1's$ in any two blocks of length $n$ differs by at  most one. There is an obvious equivalent definition for bi-infinite sequences.

If $(s_i)$ is Sturmian, then the points of the orbit of $x=0.s_1s_2 \ldots$ (binary) under $f(z)=z^2$ on the unit circle are necessarily in the same order as the points of some rigid rotation $R_{\theta}$, and vice versa. This $\theta$ is uniquely determined; it is by definition the \emph{rotation number} of $(s_i).$ Equivalently, $\theta$ is the limiting frequency of $1's$  in the sequence \cite{BS94}.

For each rational $p/q$ (modulo $1$) in lowest terms, there is a unique (up to cyclic permutation) finite word $W_{p/q}=(s_i)\in\{0,1\}^q$ such that the orbit of $x=0.\overline{s_1\ldots s_q}$ under $f(z)=z^2$ is in the same order around the circle as the points of an orbit of the rigid rotation $R_{p/q}$ (here $\overline{s_1\ldots s_q}$ denotes a recurring block). For example $W_{1/3}=001$, and $W_{2/5}=00101$.

We call $W_{p/q}$ the finite Sturmian word of rotation number $p/q$, since
the bi-infinite sequence made up of repeated copies of $W_{p/q}$ is the unique (up to shift) periodic Sturmian sequence of rotation number $p/q$. Finally we remark that there is nothing special about the symbols $1$ and $0$: identical terminology for Sturmian sequences and words may be applied if we replace $1$ and $0$ by
$\alpha$ and $\beta$ respectively.\\

We now return to the situation that $a\in \mathcal{C}_{\Gamma}$, and $z$ is a repelling fixed point of $f_a: \Lambda_{a,-} \to \Lambda_{a,-}.$ If $\hat\gamma$ is a periodic geodesic landing at $z$, it has a combinatorial rotation number $p/q$ (by Theorem \ref{lp}), and any finite word $W$ in $\alpha$ and $\beta$ which fixes $\varphi_a\circ\hat\gamma$ is Sturmian, hence (a cyclic permutation of) a power of $W_{p/q}$. By establishing and applying bounds for the eigenvalues of the Sturmian words $W_{p/q}$ in $\alpha$ and $\beta$, we 
prove our Yoccoz inequality, Theorem \ref{xxx} (see \cite{BL17}).\\

\section{Hyperbolic correspondences} \label{tra}

 We now turn to the study of the one parameter family of holomorphic correspondences defined by \eqref{fxs}.  This family is perhaps the simplest generalization of the quadratic family as a multifunction. 
 
 It will be useful  to recall some well-known facts directly related to the dynamics of $\mathbf{f}_c(z)=z^{\beta} +c$ when $\beta>1$ is a rational number.

\paragraph{Hyperbolic quadratic maps.} The notion of hyperbolicity can be given in several equivalent forms. According to the simplest one,   $f_c(z)=z^2+c$ is \emph{hyperbolic} if $f_c^n(0)$ converges to an attracting cycle (finite or infinite).

Since every finite attracting cycle attracts the orbit of a critical point, the map $f_c$ can have at most one finite attracting cycle. 
 Any quadratic map with a finite attracting cycle corresponds to a point in the interior of the Mandelbrot set $M,$  and an equivalent form of the Fatou conjecture states that this is the only possibility for a quadratic map in the interior of $M. $

On the other hand, if $c$ is in the complement of $M, $ then $\mathcal{J}_c$ is a Cantor set and $f_c$ is hyperbolic because $f_c^n(0) \to \infty.$

 The closure of attracting cycles is denoted by $\mathcal{J}_c^*.$  It turns out that \emph{$f_c$ is hyperbolic iff the basin of attraction of $\mathcal{J}_c^*$ is $\widehat{\mathbb C} \setminus \mathcal{J}_c.$  } For this reason, we call $\mathcal{J}_c^*$ the dual Julia set of $f_c.$

 This equivalent definition of hyperbolicity should be preserved in any generalization, mainly because of its intrinsic dynamical significance. 
  
 We shall use this equivalent property to define hyperbolic correspondences and centers  in the family $\mathbf{f}_c(z) = z^{\beta}+c,$  but first we need to extend the concepts of orbit, Julia set and multiplier of a cycle.

\paragraph{Cycles.}   Consider the family \eqref{fxs}.  Every sequence $(z_i)_0^{\infty}$ for which the points satisfy  $z_{i+1}\in \mathbf{f}_c(z_i)$  is a \emph{forward orbit}. A backward orbit is characterized by $z_{i+1} \in \mathbf{f}_c^{-1}(z_i).$ If $\varphi: U \to \mathbb{C}$ is an injective holomorphic map from a region $U$ of the plane such that $\varphi(z) \in \mathbf{f}_c(z),$ for every  $z$ in  $U,$ then $\varphi$ is a \emph{univalent branch} of $\mathbf{f}_c.$  By a cycle we mean any periodic forward orbit with minimal period $n.$ The quantity $$\lambda = \prod_{0}^{n-1} \varphi_i'(z_i), $$ where $\varphi_i$ is the unique univalent branch taking $z_i$ to $z_{i+1},$ is the \emph{multiplier} of the cycle. If $z=0$ then there is no univalent branch defined at $z;$ if some point of the cycle is $0,$ then by definition $\lambda=0.$

The cycle is \emph{repelling} if $|\lambda|> 1,$ and attracting if $|\lambda|< 1.$

\paragraph{Julia sets.} The Julia set of $\mathbf{f}_c$, denoted by $\mathcal{J}_c$, is the closure of the union of all repelling cycles of $\mathbf{f}_c.$ Similarly, the dual Julia set $\mathcal{J}_c^*$ is the closure of the union of all \emph{finite} attracting cycles.  The dual Julia set containing the attracting fixed point $\infty$ is denoted  by $\mathcal{J}_c^{e*}=\mathcal{J}_c^*\cup \{\infty\}.$

   \begin{figure}[H]
\centering
\includegraphics[scale=0.33]{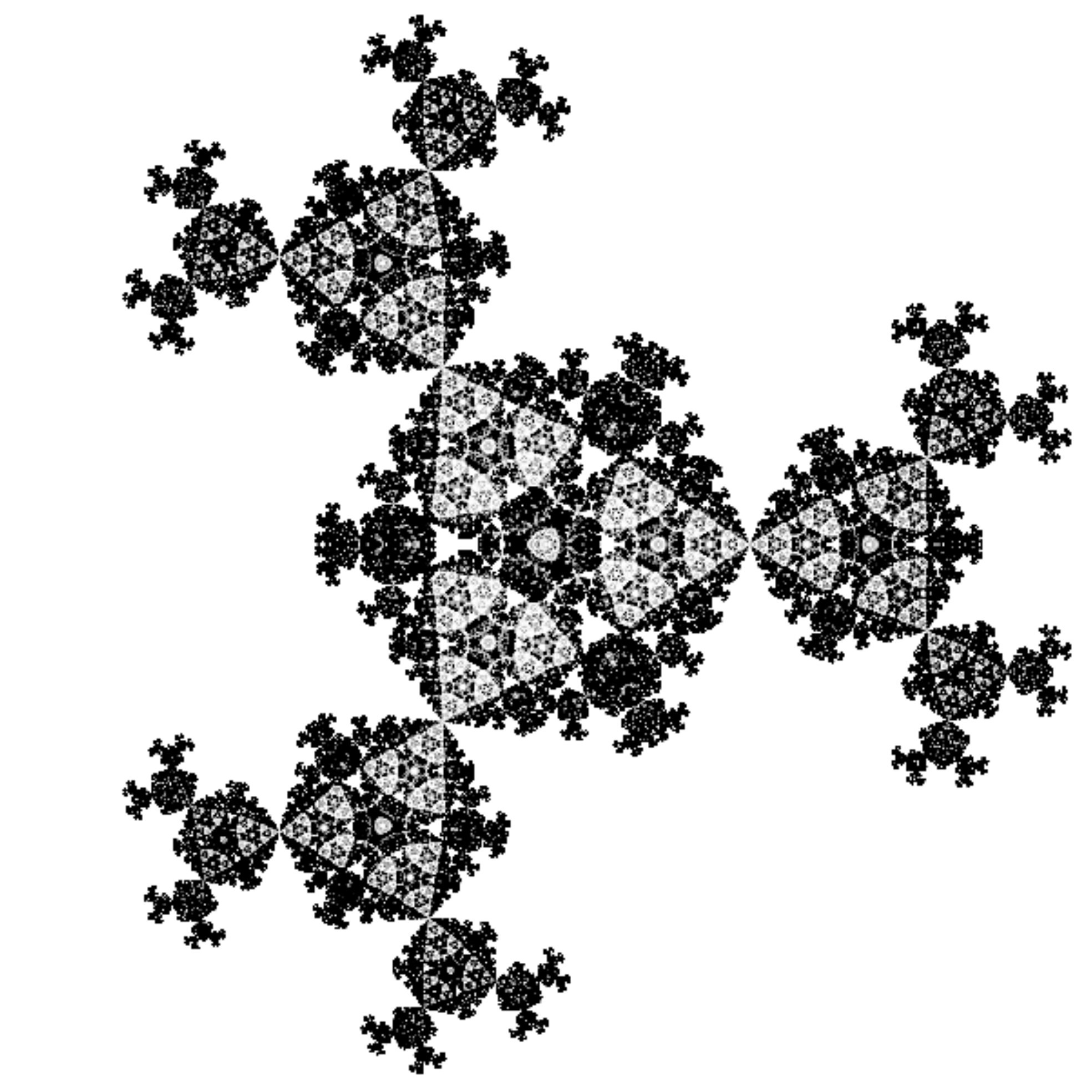}
\caption{ The Julia set of $z\mapsto z^{\frac{3}{2}} +2.$    }\label{fff}
\end{figure}

\paragraph{Filled Julia set.}

For every $c$ there is bounded disk $B$ centered at $0$ whose complement is invariant under $\mathbf{f}_c,$ and every forward orbit of a point in $\mathbb{C} \setminus B$ converges exponentially fast to $\infty.$  

We define \begin{equation} K_c = \bigcap_{n>0
} \mathbf{f}_c^{-n}(B)
\end{equation}
 as the \emph{filled Julia set} of $\mathbf{f}_c.$ A point $z$ belongs to $K_c$ iff there is at least one bounded forward orbit under $\mathbf{f}_c$ starting at $z.$ The restriction $\mathbf{f}_c|_{K_c}$ is denoted by $\mathbf{g}_c: K_c \to K_c.$

\paragraph{Hyperbolic correspondences.} 
The $\omega$-limit set of a  point $z$, denoted $\omega(z)$, consists of every $\zeta$ such that $z_{i_k} \to \zeta$ as $k\to \infty,$ for some bounded forward orbit $(z_i)$ starting at $z_0=z,$ and some subsequence $(z_{i_k}).$  We may use $\omega(z,\mathbf{f}_c)$ to make explicit the dependence on the dynamics of $\mathbf{f}_c.$

The dual Julia set is a \emph{hyperbolic attractor} for $\mathbf{g}_c$  if $\mathcal{J}_c^*$ is $\mathbf{g}_c$-forward invariant and supports an attracting conformal metric $\rho(z)|dz|,$ in the sense that

$$\sup_{z, \varphi}\| \varphi'(z) \|_{\rho} <1, $$ where the $\sup$ is taken over all $z\in \mathcal{J}_c^*$  and all univalent branches $\varphi$ of $\mathbf{f}_c$ at $z$ such that $\varphi(z) \in \mathcal{J}_c^*.$   It is implicit in this definition that $
\mathcal{J}_c^*$ does not contain the critical point, for then no univalent branch is defined at $0.$

If  $\mathcal{J}_c^*$  is a hyperbolic attractor for $\mathbf{g}_c,$ then the \emph{basin of attraction}  of $\mathcal{J}_c^{*e}$ is well defined (in other words, it contains a neighborhood of $\mathcal{J}^{*e}_c$) and consists of all  $z$ such that $\omega(z)\subset \mathcal{J}_c^{*e}.$ 

\begin{defi}[Hyperbolicity] We say that $\mathbf{f}_c$ is \emph{hyperbolic} if $\mathcal{J}_c^*$ is a hyperbolic attractor for $\mathbf{g}_c$ and the  basin of attraction of $\mathcal{J}_c^{e*}$ is  $\hat{\mathbb{C}} \setminus \mathcal{J}_c. $
 \end{defi}

\paragraph{Carpets and connectedness locus.}

 A  connected compact subset  of the plane is \emph{full} if its complement in the Riemann sphere is connected.

 A set $\Lambda\subset \mathbb{C}$  is a  hyperbolic repeller of $\mathbf{f}_c$ if (i) $\mathbf{f}_c^{-1}(\Lambda)=\Lambda;$ and (ii) $\Lambda$ supports an expanding conformal metric defined on a neighborhood of $\Lambda.$ (See \cite{SS17}). A filled Julia set $K_c$ is a \emph{Carpet} if (i) $K_c$ is connected but not full; and (ii) $K_c$ is a hyperbolic repeller. 

Intuitively, every Carpet presents holes, and by the contraction of the branches of $\mathbf{f}_c^{-1}$, every hole comes with infinitely many small copies.

We say that $K_c$ is a \emph{Cantor repeller} if $K_c$ is a hyperbolic repeller and also a Cantor set. In this case, $\mathcal{J}_c=K_c.$

 \noindent The connectedness locus $M_{\beta}$ of the family $\mathbf{f}_c$ is by definition the set of all parameters $c$ for which $K_c$ is connected.

 Another important subset of the parameter space is \begin{equation}
 M_{\beta,0} = \{c\in \mathbb{C}: 0 \in K_c\}.
 \end{equation}
Notice that both sets  generalize the definition of Mandelbrot set for the quadratic family, but if $\beta$ is not an integer, there is no reason to believe that $M_{\beta} = M_{\beta,0}.$

\begin{figure}[H] 
\centering
\includegraphics[scale=0.5]{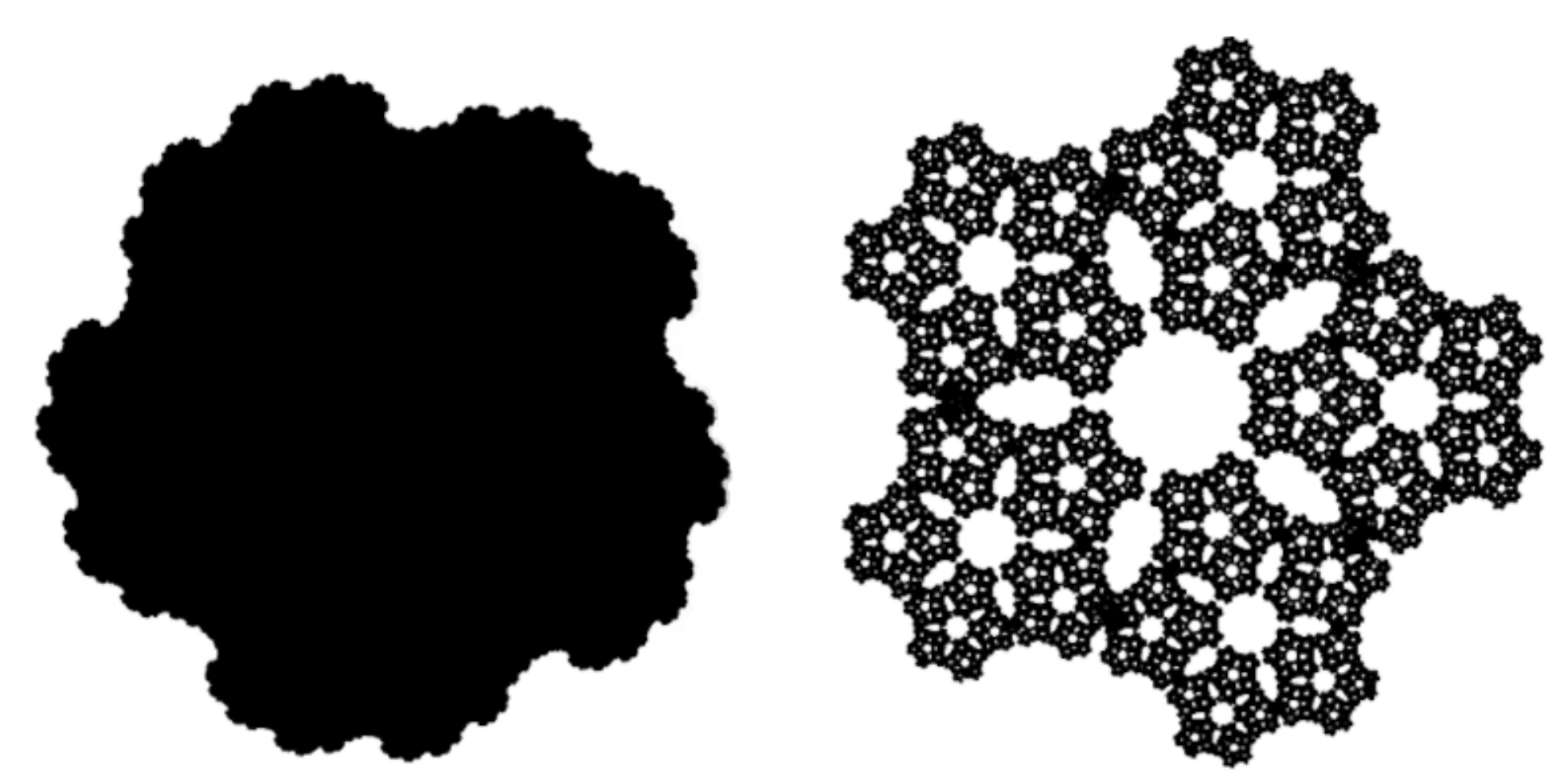} 
\caption{ Filled Julia sets in the family $\mathbf{f}_c(z)=z^{\beta} +c,$ where $\beta =5/4.$ In the  first figure (left), $K_c$ is a full compact set corresponding to  $c=3+2i.$ In the second we have a Carpet for $c= 26.$ If $|c|$ is sufficiently large, $K_c$ is a Cantor repeller.  Since the Mandelbrot set $M$ is contained in a disk of radius 2 around the origin, the fact that $c=26$ and $\mathcal{J}_c$ is still connected seems odd. However, this is one of the main features of the family $\eqref{fxs},$ and  it is experimentally clear that $M_{\beta}$  tends to cover the plane as $\beta \to 1+.$ 
}
\label{pcy}
\end{figure}

 \begin{thm} \label{pcv} If $\beta=p/q$ and  $p$ is prime, then $K_c$ is either full, a Carpet, or a Cantor repeller. 
 \end{thm}

If $p$ is prime, it is possible to show that $M_{\beta,0} \subset M_{\beta},$ in other words, $K_c$ is connected if $0 \in K_c.$ If $c$ is in $M_{\beta} - M_{\beta,0}$ then $K_c$ is a Carpet, and if $c$ is in the complement of $M_{\beta},$ then $K_c$ is a Cantor repeller.

\paragraph{Centers.}  A center is a point $c$ of the parameter space such that $$\mathbf{g}_c^n(0)=\{0\},$$ for some $n>0.$ This  definition is motivated by a well-known fact from the quadratic family, where every bounded hyperbolic component $U$ has a center   \cite{DH84, DH85}  defined as the unique point $c\in U$ for which the multiplier of the finite attracting cycle of $f_c$ is zero.  
%The center of the unbounded hyperbolic component   $\hat{\mathbb{C}} \setminus M$ is by definition $\infty.$ 

Hence, in the case of the quadratic family, the number  of bounded hyperbolic components is countably infinite, and every such component  is encoded by  a solution of $f_c^n(0)=0,$  for some $n>0.$

\paragraph{Simple centers.}   A center is called \emph{simple} if there is only one orbit of $0$ under $\mathbf{g}_c,$ and this orbit is necessarily a cycle containing $0.$ 

Let $\mathcal{S}_d = \{a\in \mathbb{C}: a^{d-1}=-1\},$ for $d>1.$ For every pair $(d,a)$ in the infinite set
 $$\bigcup_{d>1} \{d\} \times \mathcal{S}_d, $$ the point $a$ is a simple center of family of the  holomorphic correspondences $\mathbf{f}_c: z\mapsto w$ given by  $(w-c)^{2} = z^{2d}.$ Indeed, it was shown in \cite{SS17} that  the first two iterates of $0$ under $\mathbf{f}_a$ are $0 \mapsto a \mapsto a^n +a=0$ and
 $0\mapsto a \mapsto -a^{n} +a=-2a^n,$ where $-2a^{n}$ is a point in the basin of infinity of $\mathbf{f}_a.$

\paragraph{Open problems.} A fundamental program for the family $\mathbf{f}_{c}(z)=z^{\beta} +c$ is given by the following problems:
\begin{itemize}
\item[I.]  Show that every perturbation of a center corresponds to a hyperbolic correspondence;
\item[II.]  Show that the set $M_{\beta}'$ of hyperbolic parameters is indeed open and  every component of $M_{\beta}'$ is encoded by a center;
\item[III.]  Decide if the set of parameters for which $c\mapsto \mathcal{J}_c$ is continuous in the Hausdorff topology is open and dense (computer experiments seem to support this statement);

\item[IV.]  Show that every component of $\mathbb{C}\setminus \partial M_{\beta}$ is hyperbolic.

\item[V.] Classify Julia sets with zero Lebesgue measure.

\end{itemize}

The first Problem I can be solved with a generalization of the proof of Theorem \ref{jkl} (see \cite{Rigidity} for a detailed exposition);  the second is very realistic but still unresolved; the third is in many aspects  a generalization of the celebrated work of Ma\~n\'e, Sad and Sullivan \cite{Mane1983} (see also \cite{SS17} and Section \ref{dfc} for a discussion of holomorphic motions in the family \eqref{fxs}) ; and the fourth and fifth may be as difficult as the Fatou conjecture (which has been open for a century).  Indeed, the Fatou conjecture is equivalent to the following assertion \cite{Mane1983}: if $c$ is in the interior of the Mandelbrot set, then the Julia set of $f_c(z)=z^2 +c$ has zero Lebesgue measure. Theorem \ref{tyu} is perhaps the first result  towards this classification.

\begin{thm}[Hyperbolicity] \label{jkl} If $c$ is in the complement of $M_{\beta,0},$ or $c$ is sufficiently close to a simple center, then $\mathbf{f}_c$ is hyperbolic. 
\end{thm}

\subsection{Holomorphic motions} \label{dfc}

 Quasiconformal deformations of Julia sets in the family $\mathbf{f}_c$ can be explained  by the theory of branched holomorphic motions  introduced by Lyubich and Dujardin \cite{Lyubich2015} for polynomial automorphisms of $\mathbb{C}^2.$ For more details, see \cite{SS17}.

First, let us recall some classical facts about holomorphic motions.

 Let $\Lambda \subset \mathbb{C}^n$ and $U\subset \mathbb{C}$ be an open set. A family of injections $h_c: \Lambda \to \mathbb{C}^n$ is a \emph{holomorphic motion} with base point $a\in U$ if (i) $h_a$ is the identity, and (ii) $c\mapsto h_c(z)$ is holomorphic on $U$, for every $z$ fixed in $\Lambda.$

 \paragraph{Branched holomorphic motions.} 
 
 Let $\Lambda$ and $U$ be  subsets of $\mathbb{C}$ and suppose $U$ open and nonempty. A branched holomorphic motion with base point $a\in U$ is a multifunction $\mathbf{h}: U \times \Lambda \to  \mathbb{C}$ with the following properties: (i) $\mathbf{h}(a,z)=\{z\},$ for every $z\in \Lambda.$ In other words, $\mathbf{h}_a=\mathbf{h}(a, \cdot)$ is the identity; and
(ii) there is a family $\mathcal{F}$ of holomorphic maps $f:U \to \mathbb{C}$ such that $$\bigcup_{z\in \Lambda} G_z(\mathbf{h}) = \bigcup_{f\in \mathcal{F}} G(f), $$
where $G(f)=\{(z, fz);\,  z\in U\}$ is the graph of $f$ and $G_{z}(\mathbf{h})$ is the graph of $c\mapsto \mathbf{h}_c(z).$

The key difference in the definitions of branched and (non-branched) holomorphic motion is that bifurcations are allowed in the branched family, so that $\mathbf{h}_c(z)$ is a set instead of a single point.

 \subsection{Solenoidal Julia sets.} Recently, Siqueira and Smania have presented another way of interpreting branched holomorphic motions on the plane as projections  of (non-branched) holomorphic motions on $\mathbb{C}^{2}.$ The method is general and applies to every hyperbolic Julia set \cite{SS17}, but we shall restrict to bifurcations near $c=0.$

 There is a family of holomorphic maps $f_c: U_0 \to V_0$ such that $U_0$ and $V_0$ are open subsets of $\mathbb{C}^2,$ the closure of $U_0$ is contained in $V_0,$ and the maximal invariant set
 $$ S_c = \bigcap_{n=1}^{\infty} f_c^{-n}(V_0) $$

  \noindent is the closure of periodic points of $f_c.$ (All periodic points are repelling in a certain generalized sense, see \cite{SS17}). This description holds for every $c$ in a neighborhood of zero. The dynamics of $\mathbf{f}_c$ on $\mathcal{J}_c$ is a topological factor of $f_c: S_c \to S_c,$ in the sense that $\pi(S_c)=\mathcal{J}_c$ and $\pi$ sends two points in $S_c$ related by $f_c$ to two points in $\mathcal{J}_c$ related by $\mathbf{f}_c$:  $\pi f_c(x)$ is an image of $\pi(x)$ under $\mathbf{f}_c,$ for every $x\in S_c.$

 Let $\pi_c: S_c \to \mathcal{J}_c$ denote the projection $(z,w) \mapsto z.$

 \begin{thm}[Holomorphic motions] \label{bfl} There is a holomorphic motion $h_c:S_0 \to \mathbb{C}^2$ with base point $c=0$ such that 
 
 \begin{enumerate}
 
 \item $h_c(S_0) =S_c$ and $h_c$ is a conjugacy (homeomorphism) from $f_0:S_0 \to S_0$ to $f_c:S_c \to S_c.$
 
 \item the projected motion $\mathbf{h}_c(z)= \pi_c \circ f_c \circ \pi_0^{-1}(z) $ is a branched holomorphic motion mapping $\mathcal{J}_0=\mathbb{S}^1$ to $\mathcal{J}_c=\mathbf{h}_c(\mathbb{S}^1).$

 \item $S_0$ is a solenoid, and $\mathbf{f}_c$ is hyperbolic, for every $c$ in $U.$
 
 \end{enumerate}
 
 \end{thm}

 \noindent See \cite{SS17} for the solenoidal description of $S_0$ (indeed,  $S_0$ is  the Williams-Smale solenoid for certain values of $p$ and $q$).   
 
 In Figure \ref{bcx}, the motion of $\mathcal{J}_c$ is illustrated in four steps. 
 
%%%%% Figure here
\subsection{Conformal iterated function systems}

  Dual Julia sets $\mathcal{J}_c^*$ in the family \eqref{fxs} often appear as limit sets of conformal iterated function systems (CIFS).  This phenomenon is easy to explain when $c$ is close to zero, and very convenient to motivate further generalizations.  

Indeed, using the contraction of $\mathbf{f}_c$ around $z=0$ one can prove that for every $c\neq 0$ close to zero, there is an open disk $D$ such that $D_1=\mathbf{f}_c(D)$ is another disk avoiding zero and compactly contained in $D.$ 

Since $D_1$ is simply connected, there are $q$ conformal branches $f_j: D_1 \to \mathbb{C}$ such that $\mathbf{f}_c(z)=\{f_j(z) \}_j,$ for every $z\in D_1.$ Moreover, the images $f_j(D_1)$ are disjoint disks. It follows that $$\mathbf{f}_c(D_1) \subset \mathbf{f}_c(D)=D_1; $$  and the family of maps $f_j:  D_1 \to D_1$ is a CIFS.  The limit set of this CIFS is   $\Lambda=\cap_n H^n(D_1),$ where $$H(A) = \bigcup_{j=1}^{q} f_j(A)$$ is the Hutchinson operator, and $A\subset D_1.$  The most important fact derived from this construction is that $\Lambda$ is the closure of attracting periodic orbits:   $\Lambda=\mathcal{J}_c^*.$

This analysis has many generalizations, including holomorphic motions and Hausdorff dimension. Theorem \ref{tyu}, for example, is stated in great generality in \cite{Siq17}.  

In \cite{Rigidity}  we give a general account establishing a rigidity result which states that $\mathcal{J}_c^*$ is finite at simple centers, but   any perturbation of $c$ yields a hyperbolic correspondence whose dual Julia set is a Cantor set. In the case of $c$ close to zero, for example,  $\mathcal{J}_c^*$ is either a Cantor set  if $c\neq 0$ (indeed, $\Lambda$ comes from a CIFS without overlaps) or a single point set  $\mathcal{J}_0^*=\{0\}.$

\bibliographystyle{spphys.bst}

\bibliography{oi}  

\Addresses

\end{document}